\documentclass[epsbox,12pt,reqno]{amsart}
\usepackage{latexsym}
\usepackage{amsmath,amssymb,latexsym,amscd} 
\usepackage{graphics}

\newtheorem{theorem}{Theorem}

\newtheorem{lemma}[theorem]{Lemma}
\theoremstyle{definition}
\newtheorem{definition}[theorem]{Definition}

\theoremstyle{remark}
\newtheorem{remark}[theorem]{Remark}
\theoremstyle{preliminaries}

\begin{document}
%
\title{A quadratic lower bound for colourful simplicial depth}
\author{Tamon Stephen}
\address{
Department of Mathematics,
Simon Fraser University,
8888 University Drive,
Burnaby, British Columbia, Canada  V5A 1S6.
}
\email{tamon@sfu.ca}
\author{Hugh Thomas}
\address{
Department of Mathematics and Statistics,
University of New Brunswick,
Fredericton, New Brunswick E3B 5A3
Canada.
}
\email{hugh@math.unb.ca}
\thanks{Both authors were supported by NSERC Discovery grants.
Additionally, T.~Stephen was supported by DFG FG-468 and the 
Dynamical Systems research focus at the University of Magdeburg.}


%
\def\bara{B\'ar\'any}
\def\mato{Matou\v sek}
\def\P{\mathbb{P}}
\def\R{\mathbb{R}}
\def\S{\mathbf{S}}
\def\Sph{\mathbb{S}}
\def\sd{\operatorname{depth}}
\def\csd{\mathbf{depth}}
\def\core{\operatorname{core}}
\def\conv{\operatorname{conv}}
\def\cone{\operatorname{cone}}
\def\zero{{\bf 0}}

\maketitle

\begin{abstract}
We show that any point in the convex hull of
each of $(d+1)$ sets of $(d+1)$ points in $\R^d$ 
is contained in at least 
$\left\lfloor {(d+2)^2}/{4}\right\rfloor$
simplices with one vertex from each set.
\end{abstract}


\section{Introduction}
Given a set $S$ of points in $\R^d$ and an additional point $p$,
the {\it simplicial depth} of $p$ with respect to $S$, denoted
$\sd_S(p)$,  is the
number of closed $d$-simplices generated from points of $S$ that contain $p$.
This can be viewed as a statistical measure of how representative
$p$ is of $S$ \cite{FR05}.
In \cite{DHST05} the authors consider configurations
of $d+1$ points in each of $d+1$ colours in $\R^d$.  They 
define the {\it colourful simplicial depth} of $p$ with respect to 
a configuration $\S$,
denoted $\csd_{\S}(p)$, as the number of $d$-simplices containing $p$
generated by sets of points from $\S$ that contain one point of
each colour.  

Given a configuration $\S = \{S_1, \ldots, S_{d+1}\}$ 
the {\it core} of the configuration is the intersection 
of the convex hulls of the individual colours, i.e. 
$\bigcap_{i=1}^{d+1} \conv(S_i)$.  Define: 
\begin{equation}\label{eq:mu}
\mu(d)= \min_{\text{configurations } \S \text{ in } \R^d, 
               ~ p \in \core(\S)} \csd_S(p)
\end{equation}
The quantity $\mu(d)$ was introduced in \cite{DHST05}.  In that paper,
it was shown that $2d \le \mu(d) \le d^2+1$, and conjectured
that $\mu(d)=d^2+1$.  In this paper we prove
\begin{theorem} \label{th:main}
$\mu(d) \ge \left\lfloor {(d+2)^2}/{4} \right\rfloor$.
\end{theorem}
In particular, this shows that $\mu(d)$ is quadratic.
The quantity $\mu(d)$ is used in bounding the depth of a monochrome
simplicial median (i.e.~point of maximum simplicial depth) for
$n$ points in $\R^d$ 
via the method of \bara\ \cite{Bar82} as described in \cite{DHST05}.  
We remark also that in optimization, 
$\mu(d)$ represents the minimum number of solutions
to the colourful linear programming feasibility problem
proposed in \cite{BO97} and discussed in \cite{DHST06}.

\section{Preliminaries} \label{se:prelim}
We consider only configurations that have a non-empty core.
Since we compute depths using {\it closed} simplices, 
degeneracies that cause $p$ to lie on
the boundary of a colourful simplex can only increase the 
colourful simplicial depth 
by allowing $p$ to lie in different simplices with disjoint interior.
Thus, since we are minimizing, we can assume 
that the core is full-dimensional and the points of
$\S$ lie in general position in $\R^d$. 

We also assume without loss of generality that 
the minimum in Equation~(\ref{eq:mu}) is attained at
the origin, $p=\zero$.
We note that if some point in $\S$ is $\zero$ then we are done
since all the $(d+1)^d$ colourful simplices using this point contain $\zero$.
Thus 
we can rescale the non-zero points of $\S$ so that they lie on the unit sphere,
$\Sph^d \subset \R^d$. 
Since the coefficients in a convex combination expressing $\zero$
can also be rescaled, this does not affect which colourful simplices
contain $\zero$.

Indeed, we observe that the colourful set $\{x_1,\ldots,x_{d+1}\}$
generates a colourful simplex containing $\zero$ exactly when the antipode
$-x_{d+1}$ of $x_{d+1}$ lies in $\cone(x_1,\ldots,x_d)$,
a pointed cone with vertex $\zero$.
Our strategy will be to understand how $\Sph^d$ can be covered
by $d$-coloured simplicial cones, that is, cones that are 
generated by $d$ points of different colours.
In this vein we can define the {\it D-depth} of a point of colour $i$ 
to be the number of $d$-coloured simplicial cones of colours
$D=\{1,\ldots,\hat i,\ldots,d+1\}$ containing the point.
We remark that the $D$-depth of any point is at least one.
This follows from the result in \cite{Bar82} that every point 
in a colourful configuration
with $\zero$ in its core is among the generators of at least one colourful 
simplex containing $\zero$.

Let $e_1, \ldots, e_d$ be the standard coordinate unit vectors in $\R^d$.
Recall that the {\it standard cross-polytope} 
is $\conv(\pm e_1, \ldots, \pm e_d)$.  We will now define a condition on
$2n$ points that means that they ``look like'' the vertices of a
standard cross-polytope, with $\pm e_i$ coloured with colour $i$.  

\begin{definition} A collection of 2 points in each of $d$ colours is
said to be in {\it deformed cross position} if the $2^d$
different $d$-coloured simplicial cones generated by the points
cover $\mathbb{R}^d$.
\end{definition}

Note that some of 
the $d$-coloured simplicial cones generated by the points in deformed 
cross position may overlap substantially (not just along boundaries).
We conclude with the following Lemma, which is proved in
Section~\ref{se:lemma}.

\begin{lemma} \label{le:}
If the colourful simplicial depth of $\zero$ is less than
$d^2+d$, then for any choice of a set $D$ of $d$ colours, there must exist a
subset of $\S$ in deformed cross position, the colours of whose vertices 
are given by $D$.
\end{lemma}

\section{Proof of Theorem~\ref{th:main}}
Assume that the colourful simplicial depth of $\zero$ is less 
than $d^2+d$, so that the lemma applies.

Choose a set of points $P_1$ in deformed cross position on the colours
$\{2,\dots,d+1\}$.  Pick a point $v$ from $\S$
with colour $1$.  Its antipode is
in at least one $\{2,\dots,d+1\}$-coloured simplicial cone generated 
by vertices of $P_1$.  The vertices of that cone together with $v$ yield
a colourful simplex containing $\zero$.  This procedure yields $d+1$ 
colourful simplices, one for each element of $\S$ with
colour 1.

Now choose a set of points $P_2$ in deformed cross position on the colours
$\{1,3,\dots,d+1\}$.  Let $v$ be a point from $\S$ 
with colour 2 which does not
appear in $P_1$.  There are $d-1$ of these.  
As before, each of these points, together with some
vertices from $P_2$, generate a colourful simplex containing $\zero$.
Since we are using vertices of colour 2 which were not used in the first
step, the colourful simplices generated at this step are distinct from
those generated at the first step.  This yields $d-1$ colourful simplices.  

Repeat this procedure, at the $i$-th step choosing points in deformed
cross position
on the colours $\{1,\dots,\hat i,\dots, d+1\}$, and then considering
those vertices of colour $i$ which have not appeared in any $P_j$ for
$j<i$.  This gives $d+1-2(i-1)$ new colourful simplices.
Hence the total number of colourful simplices produced is at least: 
$\displaystyle (d+1)+(d-1)+\dots= \left\lfloor {(d+2)^2}/4\right\rfloor$ 
as desired.  

\begin{remark}
This improves the lower bound of $2d$ from \cite{DHST05} starting at $d=4$.
\end{remark}

\begin{remark}
The authors have recently learned that \bara~ and \mato~ 
independently found a quadratic lower bound for $\mu(d)$ \cite{BM06}.
Their bound is $\mu(d) \ge \frac{1}{5}d(d+1)$.  They also give
a lower bound of $3d$ if $d>2$ which exceeds $(d+2)^2/4$ when
$d=3,4,5,6,7$.
\end{remark}

\subsection{Proof of Lemma~\ref{le:}} \label{se:lemma}
Without loss of generality, let $D=\{1,\dots,d\}$.  
Consider the $D$-depth of a 
point in $\Sph^d$. 
If every point were of $D$-depth at least $d$, then wherever the points coloured
$d+1$ are, each of 
their antipodes is in at least $d$ $D$-coloured simplicial cones,
and thus the depth of $\zero$ is at least $d^2+d$.  

Assuming the colourful simplicial depth of $\zero$ is less than $d^2+d$,
there is some point $x\in \Sph^d$ which is in no more than $d-1$ 
$D$-coloured cones.  Thus, we can choose a set of points $w_1,\dots,w_d$
such that $w_i$ is of colour $i$ and generates 
no $D$-coloured cone containing $x$.   
Let $z_1,\dots,z_d$ be the vertices of some 
$D$-coloured cone containing $x$, with $z_i$ of colour $i$.  

We claim that $P=\{z_i\}\cup\{w_i\}$ is in deformed cross position.
Let $\P^d$ be the union of $d$-coloured simplices on the set $P$.  
Consider the map $f$ which maps $\P^d$ to $\Sph^d$
by $x \rightarrow x/||x||$.
%
We want to show that this map is onto.  Suppose otherwise.  
Let $X$ be the simplex of $\P^d$ whose vertices are 
$\{z_1,\dots,z_d\}$. Let $Y$ be the union of the other simplices of $\P^d$.  
Let $Z=X\cap Y$ be the boundary of $X$.  

Let $A$ be the intersection of $\Sph^d$ with the $D$-coloured cone generated
by the $\{z_i\}$.  Let $B$ be the boundary of $A$.  

By definition, $f(X)=A$.  Thus, if $f$ is not onto, there is some point
$y \not\in A$ such that $y$ is not in the image of $f$.  
Also observe that $x \notin f(Y)$, by our choice of points $\{w_i\}$.  

Now, define a map $\pi$ which retracts $\Sph^d\setminus\{x,y\}$ onto $B$.  
Clearly, restricted to $Z$, 
$(\pi\circ f)|_Z=f|_Z$ is a 
homeomorphism, and generates the non-zero homology of $B$.  But $\pi\circ
f:Y\rightarrow B$ shows that $(\pi\circ  f)|_Z$ is null-homotopic, which is
 a contradiction.  

Thus $f$ must be onto, and our set of points is in deformed cross position.

\section{Acknowledgments}
The authors would like to thank the referees for comments which
improved the presentation of the paper.
%
%

\bibliographystyle{spdcg}

\providecommand{\bysame}{\leavevmode\hbox to3em{\hrulefill}\thinspace}
\providecommand{\MR}{\relax\ifhmode\unskip\space\fi MR }
\providecommand{\MRhref}[2]{%
  \href{http://www.ams.org/mathscinet-getitem?mr=#1}{#2}
}
\providecommand{\href}[2]{#2}

\end{document}